\documentclass[1p,12pt]{elsarticle}

\usepackage{graphicx}
\usepackage{amsmath,amsxtra,amssymb,latexsym, amscd,amsthm}

\theoremstyle{plain}
\newtheorem{thm}{Theorem}[section]

\theoremstyle{definition}

\newtheorem{rmk}[thm]{Remark}

\numberwithin{equation}{section}
\numberwithin{figure}{section}
\numberwithin{table}{section}

\newcommand{\M}{\operatorname{M}}
\newcommand{\Lo}{\operatorname{L}}

\begin{document}

\begin{frontmatter}
\title{\textbf{A new proof for the number of lozenge tilings of quartered hexagons}}

% input author, affilliation, address and support information as follows;
% the address should include the country, and does not have to include
% the street address

\author{Tri Lai\corref{cor1}\fnref{myfootnote1}}
\address{Institute for Mathematics and its Applications\\ University of Minnesota\\ Minneapolis, MN 55455}
\fntext[myfootnote1]{This research was supported in part by the Institute for Mathematics and its Applications with funds provided by the National Science Foundation.}
\cortext[cor1]{Corresponding author, email: tmlai@ima.umn.edu, tel: 612-626-8319}

%\date{
%\small Mathematics Subject Classifications: 05A15, 05C70, 05E99}

\begin{abstract}
It has been proven that the lozenge tilings of a quartered hexagon on the triangular lattice are enumerated by a simple product formula. In this paper we give a new proof for the tiling formula by using Kuo's graphical condensation. Our result generalizes a Proctor's theorem on enumeration of plane partitions contained in a ``maximal staircase".
 %Proctor proved a simple formula for the number of a class of plane partitions contained in a  ``maximal staircase". The result is equivalent to the enumeration of lozenge tilings of a hexagon with a maximal staircase removed from some of its vertices. In this paper we give a new proof of a generalization of the Proctor's result by using Kuo's graphical condensation.
%
%
  % keywords are optional
  %\bigskip\noindent \textbf{Keywords:} tilings, perfect matchings,  plane partitions, graphical condensation
\end{abstract}

\begin{keyword}
Tilings\sep perfect matchings \sep plane partitions \sep graphical condensation.
\MSC[2010] 05A15\sep05C70 \sep 05E99
\end{keyword}

\end{frontmatter}

\section{Introduction}

%A lattice divides the plane into fundamental regions. A (lattice) \textit{region} is a finite connected union of fundamental regions of that lattice. A \textit{tile} is the union of any two fundamental regions sharing an edge. A \textit{tiling} of the region $R$ is a covering of $R$ by tiles with no gaps or overlaps. Denote by $\M(R)$ the number of tilings of the region $R$.

A plane partition is a rectangular array of non-negative integers with weakly decreasing rows and columns. The number of plane partitions contained in a $b\times c$ rectangle with entries at most $a$ is given by MacMahon's formula  $\prod_{i=1}^{a}\prod_{j=1}^{b}\prod_{k=1}^{c}\frac{i+j+k-1}{i+j+k-2}$ (see \cite{Mac}).
As a variation of this,  Proctor proved a simple product formula for the number of plane partitions with entries at most $a$ which are contained in a shape with row lengths $b, b-1,\dotsc,b-c+1$ (see Corollary 4.1 in \cite{Proctor}).

\begin{figure}\centering
\includegraphics[width=5cm]{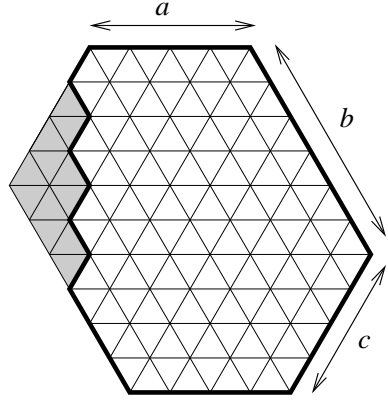}
\caption{Obtaining $P_{a,b,c}$ from the hexagon $H_{a,b,c}$ by removing a maximal staircase from the west corner.}
\label{KuoQAR5}
\end{figure}

A \textit{lozenge tiling} of a region  on the triangular lattice is a covering of the region by unit rhombi (or lozenges) so that there are no gaps or overlaps. We use notation $\Lo(R)$ for the number of lozenge tilings of a region $R$ ($\Lo(\emptyset):=1$). The plane partitions contained in a $b \times c$ rectangle with entries at most $a$  are in bijection with lozenge  tilings of  the hexagon $H_{a,b,c}$ of sides $a,b,c,a,b,c$ (in cyclic order, starting from the north side).   In the view of this we have an equivalent form of Proctor's result as follows.

\begin{thm}[Proctor \cite{Proctor}]\label{proctor} Assume that $a,b,c$ are non-negative integer so that $b\geq c$.
Let $P_{a,b,c}$ be the region obtained from the hexagon $H_{a,b,c}$ by removing the ``maximal staircase" from its east corner (see Figure \ref{KuoQAR5} for $P_{3,6,4}$).  Then
\begin{equation}
\Lo(P_{a,b,c})=\prod_{i=1}^{c}\left[\prod^{b-c+1}_{j=1}\frac{a+i+j-1}{i+j-1}\prod^{b-c+i}_{j=1}\frac{2a+i+j-1}{i+j-1}\right],
\end{equation}
where empty products are equal to $1$ by convention.
\end{thm}

One can find more variations and generalizations of the Proctor's result  in \cite{CK}. We consider next a different generalization of Theorem \ref{proctor}.

\begin{figure}\centering
\includegraphics[width=13cm]{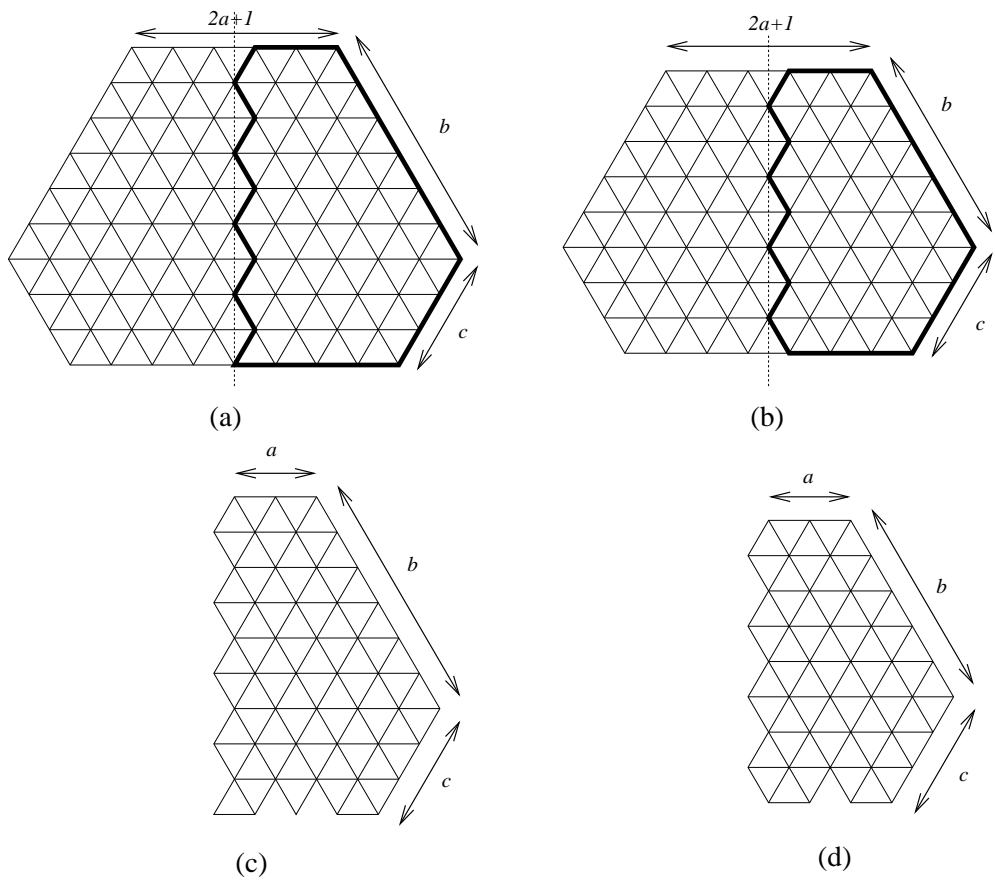}
\caption{Obtaining the region $R_{a,b,c}$ from a hexagon ((a) and (b)). Two examples of the region $R_{a,b,c}(s_1,s_2,\dotsc,s_k)$ ((c) and (d)).}
\label{KuoQAR6}
\end{figure}

\medskip

Let $R_{a,b,c}$ be the region described as in Figure \ref{KuoQAR6}. To precise, $R_{a,b,c}$ consists of all unit triangles on the right of the vertical symmetry axis of the hexagon of sides $2a+1,b,c,2a+b-c+1,c,b$ (in cyclic order, starting from the north side). Figure \ref{KuoQAR6}(a) illustrates the region $R_{2,6,3}$ and Figure \ref{KuoQAR6}(b) shows the region $R_{2,5,3}$ (see the ones restricted by the bold contours).  %If $b=0$, we let $R_{a,b,c}$ be an empty region, which has 1 ``empty tiling".
We are interested in the region $R_{a,b,c}$ with $k$ up-pointing unit triangles removed from the base ($k=\lfloor\frac{b-c+1}{2}\rfloor$). If the positions of the triangles removed  are $s_1,s_2,\dotsc,s_k$, then we denote by $R_{a,b,c}(s_1,s_2,\dotsc,s_k)$ the resulting region (see Figures \ref{KuoQAR6}(b) and (c) for $R_{2,3,6}(2,3)$ and $R_{2,5,3}(2)$, respectively). The number of lozenge tilings of the region $R_{a,b,c}(s_1,\dotsc,s_k)$ is given by the theorem stated below.

\begin{thm}\label{main} Assume $a,b,c$ are non-negative integers.
If $b-c=2k-1$ for some non-negative integer $k$, then
\begin{equation}\label{eq1}
\Lo\big(R_{a,b,c}(s_1,s_2,\dotsc,s_k)\big)=\prod_{1\leq i<j\leq k+c}\frac{s_j-s_i}{j-i}\frac{s_i+s_j-1}{i+j-1},
\end{equation}
where $s_{k+i}:= a+\frac{b-c+1}{2}+i$, for $i=1,2,\dotsc,c$.

If $b-c=2k$ for some non-negative integer $k$, then
\begin{equation}\label{eq2}
\Lo\big(R_{a,b,c}(s_1,s_2,\dotsc,s_k)\big)=\prod_{i=1}^{k+c}\frac{s_i}{2i-1}\prod_{1\leq i<j\leq k+c}\frac{s_j-s_i}{j-i}\frac{s_i+s_j}{i+j},
\end{equation}
where $s_{k+i}:= a+\frac{b-c}{2}+i$, for $i=1,2,\dotsc,c$.
\end{thm}

\begin{figure}\centering
\includegraphics[width=8cm]{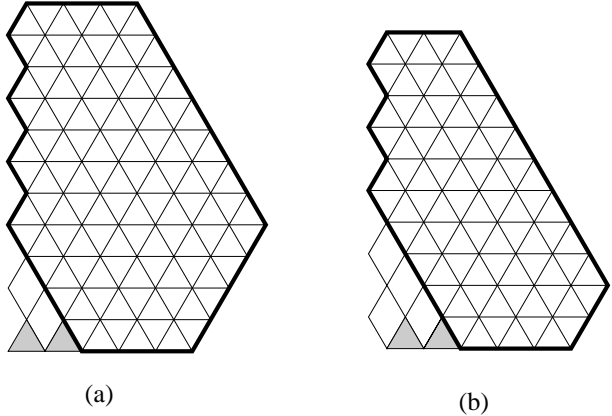}
\caption{Obtaining the region $P_{a,b,c}$ from the region $R_{a,b,c}(s_1,\dots,s_k)$.}
\label{KuoQAR7}
\end{figure}

We note that by specializing $k=b-c$ and $s_i=i$, for $i=1,2,\dotsc,k$,  the region $P_{a,b,c}$ is obtained from $R_{a,b,c}(1,2,\dotsc,k)$ by removing forced lozenges on the lower-left corner (see Figure \ref{KuoQAR7}). Thus
\[\Lo\big(R_{a,b,c}(1,2,\dotsc,k)\big)=\Lo \big(P_{a,b,c}\big),\]
and the Proctor's Theorem \ref{proctor} follows from Theorem \ref{main}.

\begin{rmk}
We enumerated the lozenge tilings of the region $R_{a,b,c}(s_1,s_2,\dotsc,s_k)$ in \cite{Tri} under the name \textit{quartered hexagon} (see Theorem 3.1, equations (3.1) and (3.2)). In particular, the region $R_{a,b,c}(s_1,s_2,\dotsc,s_k)$ is obtained from the quartered hexagon\\ $QH_{b+c,n}(s_1,s_2,\dotsc,s_k,n-c+1,n-c+2,\dotsc,n)$ ($n:=\lfloor\frac{2a+1+b+c}{2}\rfloor$) by removing several forced lozenges (see Figures 2.8 and 2.9(a) in \cite{Tri}). Thus, we still call our $R_{a,b,c}$-type regions \emph{quartered hexagons}. In \cite{Tri}, we  identified  the lozenge tilings of $R_{a,b,c}(s_1,s_2,\dotsc,s_k)$ with certain families of non-intersecting paths on $\mathbb{Z}^2$, then used Lindstr\"{o}m-Gessel-Viennot Theorem (see e.g. \cite{Lind}, Lemma 1; or  \cite{Stem}, Theorem 1.2) to turn the number of path families to the determinant of a matrix whose entries are binomial coefficients, and evaluated the determinant.
\end{rmk}

A \textit{perfect matching} of a graph $G$ is collection of edges so that each vertex of $G$ is incident to exactly one selected edge. The \textit{dual graph} $G$ of a region $R$ on the triangular lattice is the graph whose vertices are unit triangle in $R$ and whose edges connect precisely two unit triangles sharing an edge. We have a bijection between the tilings of a region $R$ and the perfect matchings of its dual graph $G$. We use notation $\M(G)$ for the number of perfect matching of a graph $G$.

In this paper, we give a new inductive proof of Theorem \ref{main} by using Proctor's Theorem \ref{proctor} as a base case. The method that we use in the proof is the \textit{graphical condensation method}  first introduced by Eric Kuo \cite{Kuo}. In particular, we will employ the following theorem in our proof.

\begin{thm}[Kuo \cite{Kuo}]\label{Kuothm}
Let $G=(V_1\cup V_2, E)$ be a planar bipartite graph, and $V_1$ and $V_2$ its vertex classes. Assume that $x,y,z,t$ are four vertices appearing on a face of $G$ in a cyclic order. Assume in addition that $a,b,c\in V_1$, $d\in V_2$, and $|V_1|=|V_2|+1$. Then
\begin{equation}\label{kuoeq}
\M(G-\{y\})\M(G-\{x,z,t\})=\M(G-\{x\})\M(G-\{y,z,t\})+\M(G-\{t\})\M(G-\{x,y,z\}).
\end{equation}
\end{thm}

\section{Proof of Theorem \ref{main}}

We only prove (\ref{eq1}), as (\ref{eq2}) can be obtained by a perfectly analogous manner.

\begin{figure}\centering
\includegraphics[width=13cm]{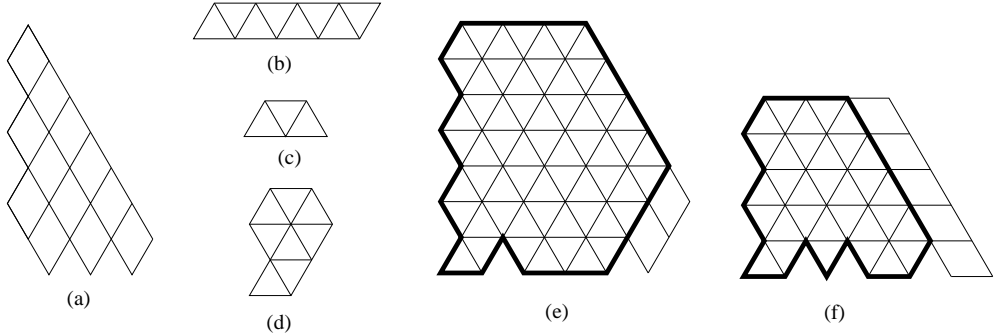}
\caption{Several special cases of $R_{a,b,c}(s_1,\dotsc,s_k)$.}
\label{KuoQAR8}
\end{figure}

It is easy to see that if $a=0$ then the region $R_{a,b,c}(s_1,s_2,\dotsc,s_k)$ has only one tiling (see Figure \ref{KuoQAR8}(a)). On the other hand, since now $\{s_1,\dotsc,s_k\}=[k]$\footnote{We use the notation $[k]$ for the set $\{1,2,\dotsc,k\}$ of all positive integers not exceed $k$.}, the right hand side of the equality (\ref{eq1}) is also equal to $1$, then (\ref{eq1}) holds for $a=0$. Moreover, if $b=0$, then $c=1$, $k=0$, and the region has the form as in  Figure \ref{KuoQAR8}(b). In this case, the region has also a unique tiling; and it is easy to verify that the right hand side of (\ref{eq1}) equals 1. Thus, we can assume in the remaining of the proof that $a,b\geq 1$.

We will prove (\ref{eq1}) by induction on $a+b$.

If $a+b\leq2$, then we have $a=b=1$. It is easy to see that there are only 2 possible shapes for the region $R_{a,b,c}$ (i.e. the region before removed triangles from the base) as in Figures \ref{KuoQAR8}(c) and (d). Then it is routine to verify (\ref{eq1}) for $a=b=1$.

For the induction step, we assume that (\ref{eq1}) is true for any region with $a+b<l$, for some $l\geq3$, then we need to show that it is also true for any region $R_{a,b,c}(s_1,\dotsc,s_k)$ with $a+b=l$.

Let $A:=\{s_1,s_2,\dotsc,s_k\}$ be a set of positive integers, we define the  operators $\Delta$ and $\bigstar$ by setting
\begin{equation}
\Delta(A):=\prod_{1\leq i<j\leq k}(s_j-s_i)
\end{equation}
and
\begin{equation}
\bigstar(A):=\prod_{1\leq i<j\leq k}(s_i+s_j-1).
\end{equation}
Then one can re-write (\ref{eq1}) in terms of the above operators as:
\begin{equation}\label{eq1a}
\Lo\big(R_{a,b,c}(s_1,\dotsc,s_k)\big)=\frac{\Delta(S)}{\Delta([k+c])}\frac{\bigstar(S)}{\bigstar([k+c])},
\end{equation}
where $S:=\{s_1,s_2,\dotsc,s_{k+c}\}$ and $[k+c]:=\{1,2,\dotsc,k+c\}$. From this stage we will work on this new form of the equality (\ref{eq1}).

\medskip

We first consider three special cases as follows:
\begin{enumerate}
\item[(i)] If $c=0$, then by considering forced lozenges as in Figure \ref{KuoQAR8}(f), we get
\begin{equation}\label{eq3}
\Lo\big(R_{a,b,0}(s_1,s_2,\dotsc,s_k)\big)=\Lo\big(R_{a-q,b-1,1}(s_1,s_2,\dotsc,s_k)\big),
\end{equation}
where $q=a+\frac{b-c+1}{2}-a_k$. Then (\ref{eq1a}) follows from the induction hypothesis for the region on the right hand side of (\ref{eq3}).

\item[(ii)] If $k=0$, then $b=c-1$; and we get the region $P_{a,b,b}$ is obtained from the region $R_{a,b,b+1}(\emptyset)$ by removing forced lozenges along its base. Thus, (\ref{eq1a}) follows from Proctor's Theorem \ref{proctor}.

\item[(iii)] Let $d:=a+\frac{b-c+1}{2}$ (so $s_{k+i}=d+i$). We consider one more special case when $a_k=d$. By removing forced lozenges again, one can transform our region into the region $R_{a,b-1,c+1}(s_1,\dotsc,s_{k-1})$, then we get again (\ref{eq1a}) by induction hypothesis for the latter region (see Figure \ref{KuoQAR8}(e)).
\end{enumerate}

\medskip

From now on, we assume that our region $R_{a,b,c}(s_1,\dotsc,s_k)$ has the two parameter $k$ and $c$ positive (so  $b=c+2k-1\geq 2$), and that $a_k < d$.

\begin{figure}\centering
\includegraphics[width=8cm]{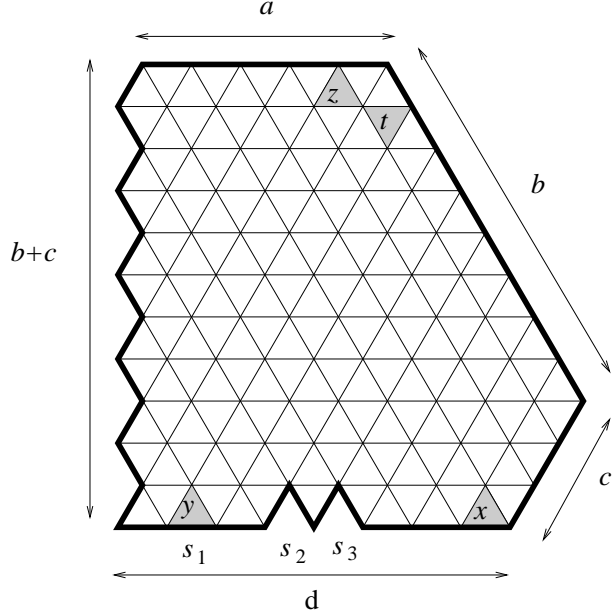}%
\caption{Region to which we apply Kuo's graphical condensation}
\label{KuoQAR4}
\end{figure}

Now we consider the region $R$ obtained from $R_{a,b,c}(s_1,\dotsc,s_k)$ by recovering the unit triangle at the position $s_1$ on the base. We now apply Kuo's Theorem \ref{Kuothm} to the dual graph $G$ of $R$, where the unit triangles corresponding to the four vertices $x,y,z,t$  are chosen as in Figure \ref{KuoQAR4} (see the shaded triangles). In particular, the triangles corresponding to $x$ and $y$ are at the positions $s_1$ and $d$ on the base; and the ones corresponding to $z,t$ are on the upper-right corner of the region.

\begin{figure}\centering
\includegraphics[width=10cm]{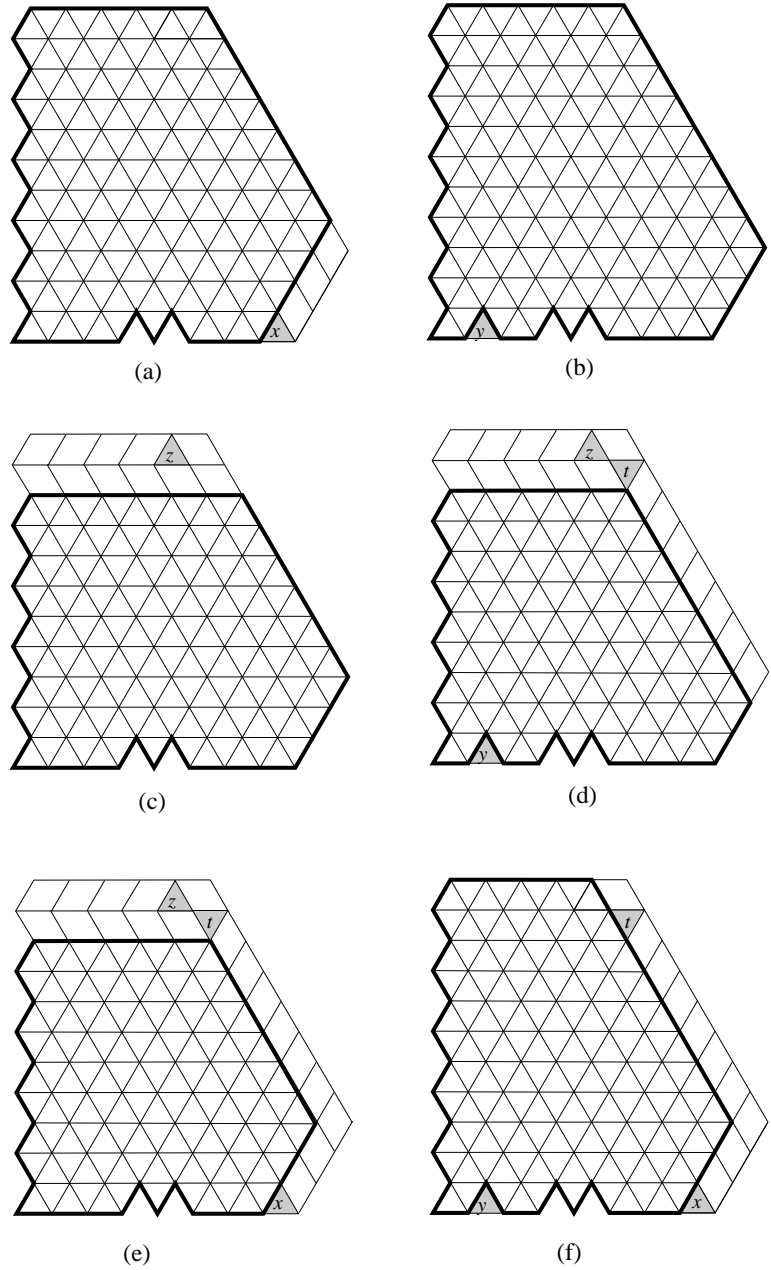}
\caption{Obtaining the recurrence for $\Lo(R_{a,b,c}(s_1,\dotsc,s_k))$.}
\label{KuoQAR9}
\end{figure}

One readily sees that the six regions that have dual graphs appearing in the equation (\ref{kuoeq}) of  Kuo's Theorem have some lozenges, which are forced into any tilings. Luckily, by removing such forced lozenges, we still get new regions of $R_{a,b,c}$-type. In particular, after removing forced lozenges from the region corresponding to $G-\{x\}$, we get the region $R_{a,b-1,c+1}(s_2,s_3,\dotsc,s_{k})$ (see the region restricted by bold contour in Figure \ref{KuoQAR9}(a)). This implies that
\begin{equation}\label{eqn1}
\M(G-\{x\})=\Lo\big(R_{a,b-1,c+1}(s_2,a_3,\dotsc,a_{k})\big).
\end{equation}
Similarly, we have five more equalities corresponding to other graphs in (\ref{kuoeq}):
\begin{equation}\label{eqn2}
\M(G-\{y\})=\Lo\big(R_{a,b,c}(s_1,s_2,\dotsc,s_{k})\big) \text{ (see Figure \ref{KuoQAR9}(b))},
\end{equation}
\begin{equation}\label{eqn3}
\M(G-\{z\})=\Lo\big(R_{a+1,b-2,c}(s_2,s_3,\dotsc,s_{k})\big) \text{ (see Figure \ref{KuoQAR9}(c))},
\end{equation}
\begin{equation}\label{eqn4}
\M(G-\{y,z,t\})=\Lo\big(R_{a,b-2,c-1}(s_1,s_2,\dotsc,s_{k})\big) \text{ (see Figure \ref{KuoQAR9}(d))},
\end{equation}
\begin{equation}\label{eqn5}
\M(G-\{x,z,t\})=\Lo\big(R_{a,b-2,c}(s_2,s_3,\dotsc,s_{k})\big) \text{ (see Figure \ref{KuoQAR9}(e))},
\end{equation}
\begin{equation}\label{eqn6}
\M(G-\{x,y,t\})=\Lo\big(R_{a-1,b,c}(s_1,s_2,\dotsc,s_{k})\big) \text{ (see Figure \ref{KuoQAR9}(f))}.
\end{equation}
Plugging the above six equalities (\ref{eqn1}) -- (\ref{eqn6}) in (\ref{kuoeq}), we have the following recurrence
\begin{align}\label{recu1}
\Lo\big(R_{a,b,c}&(s_1,s_2,\dotsc,s_{k})\big)\Lo\big(R_{a,b-2,c}(s_2,s_3,\dotsc,s_{k})\big)\notag\\&=\Lo\big(R_{a,b-1,c+1}(s_2,s_3,\dotsc,a_{k})\big)\Lo\big(R_{a,b-2,c-1}(s_1,s_2,\dotsc,s_{k})\big)\notag\\&\quad +\Lo\big(R_{a+1,b-2,c}(s_2,s_3,\dotsc,s_{k})\big)\Lo\big(R_{a-1,b,c}(s_1,s_2,\dotsc,s_{k})\big).
\end{align}
The five regions other than $R_{a,b,c}(s_1,s_2,\dotsc,s_{k})$ in the above recurrence (\ref{recu1}) have their $(a+b)$-parameter less than $l$. Therefore, by induction hypothesis, we get
\begin{equation}\label{eqm1}
\Lo\big(R_{a,b-1,c+1}(s_2,\dotsc,s_k)\big)=\frac{\Delta(S\cup\{d\}-\{s_1\})}{\Delta([k+c])}\frac{\bigstar(S\cup\{d\}-\{s_1\})}{\bigstar([k+c])},
\end{equation}
\begin{equation}\label{eqm2}
\Lo\big(R_{a,b-2,c}(s_1,\dotsc,s_k)\big)=\frac{\Delta(S-\{s_{k+c}\})}{\Delta([k+c-1])}\frac{\bigstar(S-\{s_{k+c}\})}{\bigstar([k+c-1])},
\end{equation}
\begin{equation}\label{eqm3}
\Lo\big(R_{a+1,b-2,c}(s_2,\dotsc,s_k)\big)=\frac{\Delta(S-\{s_1\})}{\Delta([k+c-1])}\frac{\bigstar(S-\{s_1\})}{\bigstar([k+c-1])},
\end{equation}
\begin{equation}\label{eqm4}
\Lo\big(R_{a-1,b,c}(s_1,\dotsc,s_k)\big)=\frac{\Delta(S\cup\{d\}-\{s_{k+c}\})}{\Delta([k+c])}\frac{\bigstar(S\cup\{d\}-\{s_{k+c}\})}{\bigstar([k+c])},
\end{equation}
\begin{equation}\label{eqm5}
\Lo\big(R_{a,b-2,c}(s_2,\dotsc,s_k)\big)=\frac{\Delta(S\cup\{d\}-\{s_1,s_{k+c}\})}{\Delta([k+c-1])}\frac{\bigstar(S\cup\{d\}-\{s_1,s_{k+c}\})}{\bigstar([k+c-1])}.
\end{equation}
By the above five equalities (\ref{eqm1}) -- (\ref{eqm5}) and the recurrence (\ref{recu1}), we only need to show that
\begin{align}\label{eqfinal}
1&=
\frac{\Delta(S\cup\{d\}-\{s_1\})\Delta(S-\{s_{k+c}\})}{\Delta(S)\Delta(S\cup\{d\}-\{s_1,s_{k+c}\})}\frac{\bigstar(S\cup\{d\}-\{s_1\})\bigstar(S-\{s_{k+c}\})}{\bigstar(S)\bigstar(S\cup\{d\}-\{s_1,s_{k+c}\})}\notag\\
&\quad\quad+\frac{\Delta(S\cup\{d\}-\{s_{k+c}\})\Delta(S-\{s_{1}\})}{\Delta(S)\Delta(S\cup\{d\}-\{s_1,s_{k+c}\})}\frac{\bigstar(S\cup\{d\}-\{s_{k+c}\})\bigstar(S-\{s_{1}\})}{\bigstar(S)\bigstar(S\cup\{d\}-\{s_1,s_{k+c}\})},
\end{align}
and (\ref{eq1a}) follows.

\medskip

First, we want to simplify the first ratio in the first term on the right-hand side of (\ref{eqfinal}). We re-write it as
\begin{equation}
\frac{\Delta(S\cup\{d\}-\{s_1\})\Delta(S-\{s_{k+c}\})}{\Delta(S)\Delta(S\cup\{d\}-\{s_1,s_{k+c}\})}=\frac{\frac{\Delta(S\cup\{d\}-\{s_1\})}{\Delta(S)} \frac{\Delta(S-\{s_{k+c}\})}{\Delta(S)}}{\frac{\Delta(S\cup\{d\}-\{s_1,s_{k+c}\})}{\Delta(S)}}
\end{equation}
The ratio $\frac{\Delta(S\cup\{d\}-\{s_1\})}{\Delta(S)}$ has its numerator and denominator almost the same, except for some terms involving $s_1$ or $d$. Canceling out all common terms of the numerator and the denominator, we have
\[\dfrac{\Delta(S\cup\{d\}-\{s_1\})}{\Delta(S)}=\frac{\prod_{i=2}^{k}(d-s_i)\prod_{i=1}^{c}(s_{k+i}-d)}{\prod_{i=2}^{k+c}(s_i-s_1)}.\]
Similarly, we get
\[\frac{\Delta(S-\{s_{k+c}\})}{\Delta(S)}=\frac{1}{\prod_{i=1}^{k+c-1}(s_{k+c}-s_i)}\]
and
\[\frac{\Delta(S\cup\{d\}-\{s_1,s_{k+c}\})}{\Delta(S)}=\frac{\prod_{i=2}^{k}(d-s_i)\prod_{i=1}^{c-1}(s_{k+i}-d)}{\prod_{i=2}^{k+c}(s_i-a_1)\prod_{i=2}^{k+c-1}(s_{k+c}-s_i)}.\]
Thus, we obtain
\begin{equation}\label{eqsimple1}
\frac{\Delta(S\cup\{d\}-\{s_1\})\Delta(S-\{s_{k+c}\})}{\Delta(S)\Delta(S\cup\{d\}-\{s_1,s_{k+c}\})}=\frac{\frac{\Delta(S\cup\{d\}-\{s_1\})}{\Delta(S)} \frac{\Delta(S-\{s_{k+c}\})}{\Delta(S)}}{\frac{\Delta(S\cup\{d\}-\{s_1,s_{k+c}\})}{\Delta(S)}}=\frac{s_{k+c}-d}{s_{k+c}-s_1}.
\end{equation}
By the same trick, we can simply the first ratio in the second term on the right-hand side of (\ref{eqfinal}) as
\begin{equation}\label{eqsimple2}
\frac{\Delta(S\cup\{d\}-\{s_{k+c}\})\Delta(S-\{s_{1}\})}{\Delta(S)\Delta(S\cup\{d\}-\{s_1,s_{k+c}\})}=\frac{\frac{\Delta(S\cup\{d\}-\{s_{k+c}\})}{\Delta(S)}\frac{\Delta(S-\{s_{1}\})}{\Delta(S)}}{\frac{\Delta(S\cup\{d\}-\{s_1,s_{k+c}\})}{\Delta(S)}}=\frac{d-s_1}{s_{k+c}-s_1}.
\end{equation}

\medskip

Next, we simply the second ratio in each term on the right-hand side of (\ref{eq1a}). By replacing the operator $\Delta$ by the operator $\bigstar$, the whole simplifying-process works in the same way  with each term $(s_j-s_i)$ replaced by $(s_i+s_j-1)$. Thus, we get
\begin{equation}\label{eqsimple3}
\frac{\bigstar(S\cup\{d\}-\{s_1\})\bigstar(S-\{s_{k+c}\})}{\bigstar(S)\bigstar(S\cup\{d\}-\{s_1,s_{k+c}\})}=\frac{s_{k+c}+d-1}{s_{k+c}+s_1-1}
\end{equation}
and
\begin{equation}\label{eqsimple4}
\frac{\bigstar(S\cup\{d\}-\{s_{k+c}\})\bigstar(S-\{s_{1}\})}{\bigstar(S)\bigstar(S\cup\{d\}-\{s_1,s_{k+c}\})}=\frac{d+s_1-1}{s_{k+c}+s_1-1}.
\end{equation}
Finally, by (\ref{eqsimple1})--(\ref{eqsimple4}), we can simplify the equation (\ref{eqfinal}) to
\begin{equation}
1=\frac{s_{k+c}-d}{s_{k+c}-s_1}\frac{s_{k+c}+d-1}{s_{k+c}+s_1-1}+\frac{d-s_1}{s_{k+c}-s_1}\frac{d+s_1-1}{s_{k+c}+s_1-1},
\end{equation}
which is obviously true with $s_{k+c}=d+c$. This means that (\ref{eqfinal}) holds, and so does (\ref{eq1}).

\section{Acknowledgements}
The author wants to thank Christian Krattenthaler for helpful comments, and Ranjan Rohatgi for useful discussion.

\end{document}